\documentclass[12pt]{amsart}

\usepackage{comments}
\newComments\DL{DL}{red}
\newComments\AL{AL}{red}
\newComments\BJ{BJ}{blue}
\newComments\FW{FW}{green}
\usepackage{soul}

\newComments\Ol{Ol}{green}

\def\cal{\relax}
\usepackage{DLdef}

\newcommand {\fbr}{{\mathfrak{br}}}

\newcommand{\un}{\underline{N}}
\newcommand{\uN}{\underline{N}}

\newcommand {\wg}{{\wedge}}

\begin{document}

\title[Deformations of $\mathfrak{o}(5)$ in characteristics $3$ and $2$]
{Deformations of the Lie algebra $\mathfrak{o}(5)$ in
characteristics $3$ and $2$}

\author{Sofiane Bouarroudj${}^1$, Alexei
Lebedev${}^2$, Friedrich Wagemann${}^3$}

\address{${}^1$Department of Mathematics, United Arab
Emirates University, Al
Ain, PO. Box: 17551; Bouarroudj.sofiane@uaeu.ac.ae\\
${}^2$Nizhegorodskij Univ. 
RU-603950 Russia, Nizhny Novgorod, pr. Gagarina 23;
yorool@mail.ru\\
${}^3$Laboratoire de Math\'ematiques Jean Leray, UMR 6629 du CNRS,
Universit\'e de Nantes, 2 rue de la Houssini\`ere, 44322 France;
wagemann@math.univ-nantes.fr }



\begin{abstract} The finite dimensional simple modular Lie algebras
with Cartan matrix cannot be deformed if the characteristic $p$ of
the ground field is equal to $0$ or greater than $3$. If $p=3$, the
orthogonal Lie algebra $\mathfrak{o}(5)$ is one of the two simple
modular Lie algebras with Cartan matrix that have deformations (the
Brown algebras $\mathfrak{br}(2; \alpha)$ are among these
$10$-dimensional deforms and hence are not counted separately); the
$29$-dimensional Brown algebra $\mathfrak{br}(3)$ is the only other
simple Lie algebra with Cartan matrix that has deformations.
Kostrikin and Kuznetsov described the orbits (isomorphism classes)
under the action of the group $O(5)$ of automorphisms of
$\mathfrak{o}(5)$ on the space
$H^2(\mathfrak{o}(5);\mathfrak{o}(5))$ and produced representatives
of the isomorphism classes. Here we explicitly describe global
deforms of $\mathfrak{o}(5)$ and of the simple analog of this
orthogonal Lie algebra in characteristic $2$.
\end{abstract}

\thanks{We are thankful to D.~Leites who raised the problem,
to P.~Grozman for his wonderful package \texttt{SuperLie}, and to
both of them for help. The first author was partially supported by
Universit\'e de Nantes.}


\maketitle

\markboth{\itshape Sofiane Bouarroudj\textup{,} Alexei
Lebedev\textup{,} Friedrich Wagemann}{{\itshape On deformations of
$\mathfrak{o}(5)$ in characteristics $3$ and $2$}}

\thispagestyle{empty}

\section{Introduction}\texttt{}

In what follows fields of positive characteristic $p$ are supposed
to be algebraically closed. The algebraic closedness is needed since
computing cohomology we assume that all algebraic equations have
solutions. This paper can be considered as an elucidation of a very
interesting but too short paper \cite{KK} and a first step in
getting results analogous to those of  \cite{KK} for $p=2$.

\ssec{Setting} It is well known (\cite{Fu, FF2}) that the finite
dimensional simple Lie algebras over $\Cee$ are {\it rigid} meaning
that they do not have non-trivial (even infinitesimal) deformations.
Over fields of characteristic $p>3$, the same is true not for all
finite dimensional simple Lie algebras, but is true for those which
possess a Cartan matrix (for a precise definition of this notion,
see \cite{BGL}), see \cite{Ru}.

For $p=3$, Rudakov and Kostrikin (\cite{Kos, Ru} also cited in
\cite{WK} as an example in the classification of simple modular Lie
algebras with Cartan matrices) introduced a parametric family of Lie
algebras containing $\fo(5)$; the generic member of this family,
denoted $L(\eps)$, has Cartan matrix
\[
 \mat{
 2 & -1 \\
 -2 & 1-\eps
  },
\text{\ where $\eps\neq 0$.}
\]
The Lie algebra $L(\eps)$ is isomorphic to $\fo(5)\simeq \fsp(4)$
for $\eps=-1$, to the Brown algebra $\fbr(2;\alpha)$ for
$\eps=1+\frac{1}{\alpha}$, where $\eps=1$ is also possible, see
\cite{WK, BGL}. Actually, the parameter $\alpha$ is more convenient
than $\eps$; the value $\eps=0$ is excluded because $L(0)\simeq
\fpsl(3)$ is of dimension $8$, so differs drastically from the other
members of the parametric family.

It soon became clear that the deformation of $\fo(5)$ depends on
more than one parameter:

Kostrikin and Kuznetsov \cite{KK} considered a certain $3$-parameter
family of Lie algebras $L(\eps, \delta,\rho)$, which includes the
$1$-parameter family $L(\eps)=L(\eps,0,0)$. The family $L(\eps,
\delta,\rho)$ was explicitly constructed in \cite{Kos} (we reproduce
it in sec. 3.1), but the mechanism producing it remained unclear.
There is, however, a construction of the same family with clear
origin.

Earlier, using the lucid explicit description of the 3-parameter
family $\Tee(a,b,c)$ of irreducible 3-dimensional $\fsl(2)$-modules
in characteristic 3 due to Rudakov and Shafarevich \cite{RSh}
Rudakov (cited in \cite{Kos}) constructed the Cartan prolong
$L(a,b,c):=(\Tee(a,b,c), \fgl(2))_*$ which is a 3-parameter family
of deformations of $\mathfrak{o}(5)$. This construction is lucid;
for explicit expressions, see \cite{GL}; both the depth and height 1
grading of $\mathfrak{o}(5)$ in its realization as
$\mathfrak{sp}(4)$ are preserved in $L(a,b,c)$. Although nobody
bothered to express $(\eps, \delta,\rho)$ in terms of $(a,b,c)$ or
the other way round, the families of deforms are isomorphic, as
noted already in \cite{Kos}. (Note that Grozman and Leites \cite{GL}
found the exceptional values of parameters $(a,b,c)$ for which the
Cartan prolong $(\Tee(a,b,c), \fgl(2))_{*,
N}=\mathop{\oplus}\limits_{i=-1}^h\fg_i$ is of height $h>1$ and
simple; the Lie algebras $L(a,b,c)$ (that exist for all values of
parameters $(a,b,c)$ for which $\Tee(a,b,c)$ is an irreducible
$\fgl(2)$-module) are {\em partial} prolongs of height $h=1$; for
the definition of various prolongs, as well as of a shearing
parameter $\uN$ for $p>0$, see \cite{Shch}.)

Kostrikin (\cite{Kos}) proved that (for a natural explanation of the
condition in terms of reflections similar to elements of the Weyl
group, see \cite{BGL})
\begin{equation}\label{Lepsiso}L(\eps)\simeq
L(\eps')\Longleftrightarrow \eps\eps'=1\quad\text{(for $\eps\neq
\eps'$)}
\end{equation} and stated the isomorphy problem of the algebras
$L(a,b,c)$ for various values of parameters.

Kostrikin and Kuznetsov (\cite{KK}) were the first to find out
that\footnote{For an explicit form of these cocycles, see
\cite{BGL}, where this was rediscovered.} $\dim H^2(L(\eps);
L(\eps))=5$, so the family of deformations of $\mathfrak{o}(5)$
depends on at most 5 parameters. For an explanation of the meaning
of the words \lq\lq at most" here, see sec. \ref{alapp}.

Kostrikin and Kuznetsov  described the orbits under the action of
the group $O(5)$ of automorphisms of $\mathfrak{o}(5)$ (see
\cite{FG}) in the variety of Lie algebras containing
$L(-1)\simeq\mathfrak{o}(5)$, which means that Kostrikin and
Kuznetsov listed all isomorphy classes of the members of the
5-parameter family of deforms of $\mathfrak{o}(5)$. The answer is as
follows:
\begin{equation}\label{answ}
\begin{minipage}[l]{14cm}
\lq\lq For $p=3$, the Lie algebra $L(-1, -1)$ and the Lie algebras
$L(\eps)=L(1+\frac{1}{\alpha})$ with the cases $\eps=1$ included and
$\eps= 0$ excluded,  represent all the isomorphy classes of simple
10-dimensional Lie algebras deforming $\mathfrak{o}(5)$, minding
(\ref{Lepsiso})".
\end{minipage}
\end{equation}
(If $p$ were equal to 0, the radius $r$ of the sphere in the
identity representation of $O(5)$ in the 5-dimensional space
$H^2(\fo(5); \fo(5))$ would have been a natural parameter, the extra
case of $r=0$ might have occurred for $p>0$. The above answer from
\cite{KK} resembles this count.)

 Here we answer the following
natural questions arising after reading the above cited papers:

(Q1) is there a basis in the space
$H^2(\mathfrak{o}(5);\mathfrak{o}(5))$ consisting of cocycles each
of which determines a global deformation linear in the parameter?

(Q2) is there a $5$-parameter family of ($10$-dimensional simple)
Lie algebras (over a field of characteristic $3$) which includes the
families from (Q1) and which corresponds to an \so{arbitrary} linear
combination of $5$ basis cocycles from (Q1)?

The answer is ``yes'' to both questions, and we will construct the
family via obstruction theory using Massey brackets in \S3.
(Actually, the question (Q1) is already answered in affirmative in
\cite{KK} but the arguments based on algebraic geometry are indirect
and non-constructive, whereas here we prove this directly giving an
additional verification.)

\subsection{$p=2$} Are there analogs of the above results for $p=2$?
Leites told us that in the early 1970s, he suggested to divide the
last row of the standard Cartan matrix of $\fo(2n+1)$ by 2, thus
making it possible to retain simplicity for $p=2$. The algebra
$\fo(2n+1)$ itself does not, however, possess a Cartan matrix nor is
it simple if  $p=2$; it is its derived $\fo^{(1)}(2n+1)$ that does
and is, see \cite{BGL}; and it is $\fo^{(1)}(2n+1)$ that looked as a
new series of simple Lie algebras in \cite{WK} but actually was (the
derived of) the old and well-known  $\fo(2n+1)$ with
non-conventionally normalized Cartan matrix. (At the time \cite{WK}
was written, the term \lq\lq Lie algebra $\fg(A)$ with Cartan matrix
$A$" was not properly defined and was applied not only to Lie
algebras of the form $\fg(A)$ but also to their subquotients and
algebras of derivations which have no Cartan matrix, cf.
\cite{BGL}.) The infinitesimal deformations of $\fo(2n+1)$ and
$\fo^{(1)}(2n+1)$ are calculated for small values of $n$ in
\cite{BGL4}; here we describe how the infinitesimal deformation of
$\fo^{(1)}(5)$ may be integrated to a multiparameter family of Lie
algebras.

\subsection{Open problems} (a) Obtain the analog of the result of \cite{KK}
to interpret our result for $p=2$: Describe the orbits of
$Sp(4)\simeq O(5)=\faut(\fo^{(1)}(5))$ in the 4-dimensional space
$H^2(\fo^{(1)}(5); \fo^{(1)}(5))$.

(b) Explicitly describe the non-isomorphic deforms as Lie algebras
preserving geometric objects (a tensor or a distribution) both for
$p=2$ and $p=3$.

(c) In \cite{BGL4}, the infinitesimal deformations of $\fbr(3)$ are
described. Study their integrability, isomorphy classes of the
deforms, and their interpretations as in (b).

\section{Deformation and cohomology}\texttt{}

\subsection{In characteristic different from $2$}

Let $k$ be a field of any characteristic $p\neq 2$, and
$\mathfrak{g}$ a Lie algebra over $k$. A basic reference for
questions about the cohomology of Lie algebras, especially in
relation to their deformation theory, is the book by Fuchs
\cite{Fu}.

A {\it multiparameter deformation} of $\mathfrak{g}$, or {\it
multiparameter family} of Lie algebras containing $\mathfrak{g}$ as
a special member, is a Lie algebra $\mathfrak{g}_{t}$, where
$t=(t_1,\ldots,t_r)$, given by a Lie algebra structure on the tensor
product $\mathfrak{g}\otimes_k k[[t]]$ such that the Lie algebra
$\mathfrak{g}_{0}$, i.e., the one obtained when we set $t=0$, is
isomorphic to $\mathfrak{g}$ and such that $t_1,\ldots,t_r$ are
scalars with respect to the deformed bracket. {\it A posteriori} we
see that in this paper we can confine ourselves to polynomials
instead of formal power series in $t$.

The bracket in the deformed Lie algebra is of the form
\[
\begin{split}
 [x,y]_{t_1,\ldots,t_r}&=c^0(x,y)+t_1c_1^1(x,y)+\ldots+t_r c_r^1(x,y)+{}\\
 &{}+t_1^2\alpha_{1,1}(x,y)+t_1t_2\alpha_{1,2}(x,y)\ldots+t_r^2
\alpha_{r,r}(x,y) +\ldots
\end{split}
\]
for any $x,y\in\mathfrak{g}$, where $c^0(x,y):=[x,y]$ is just the
bracket of $x$ and $y$ in $\mathfrak{g}$. By linearity, it suffices
to specify the deformed bracket of elements in $\mathfrak{g}$. The
first degree conditions say that the maps
$c_i^1\colon\mathfrak{g}\otimes_k\mathfrak{g}\tto\mathfrak{g}$ must
be anti-symmetric and $2$-cocycles (with coefficients in the adjoint
module), i.e.,  for all $i=1,\ldots,r$, we have
\[
\begin{split}
 dc_i^1(x,y,z)&:=c_i^1([x,y],z)+c_i^1([y,z],x)+c_i^1([z,x],y)-{}\\
 &{}\hphantom{:}-[x,c_i^1(y,z)]-[y,c_i^1(z,x)]-[z,c_i^1(x,y)]=0.
\end{split}
\]

Two (formal) $1$-parameter deformations $\mathfrak{g}_t$ and
$\tilde{\mathfrak{g}}_t$ given by the collections $c=(c^1, c^2,
\dots)$ and $\tilde{c}=(\tilde{c}^1, \tilde{c}^2, \dots)$, where
$c^i$ and $\tilde{c}^i$ are coefficients of $t^i$, lead to
equivalent deforms (results of deformations) (i.e.,
$\mathfrak{g}_t$ and $\tilde{\mathfrak{g}}_t$ are isomorphic as
Lie algebras by an isomorphism of the form $\tau(x)={\rm
id}_{\mathfrak{g}}(x)+\mathop{\sum}\limits_{i\geq 1}\tau_i(x)t^i$
for any $x\in\mathfrak{g}$) if and only if $\tau$ links $c$ and
$\tilde{c}$ by the following formulae (for all $n>0$):
\[
 \sum_{i+j=n}\tau_i(\tilde{c}^j(x,y))=\sum_{i+j+k=n}c^i(\tau_j(x), \tau_k(y)).
\]
For the first (i.e., infinitesimal) terms, this means that two
$1$-parameter deformations are {\it infinitesimally equivalent}
(i.e., $\tau={\rm id}+t\tau_1$ and one reasons modulo $t^2$) if
and only if their $2$-cocycles differ by a coboundary. This
coboundary is nothing else than $\tau_1$. A similar statement is
true for multiparameter deformations. In particular, if two
multiparameter deformations are infinitesimally equivalent, then
the corresponding infinitesimal cocycles are linearly dependent up
to coboundaries.

For the sake of brevity, we shall recall properties of deformations
using only $1$-parameter deformations; to generalize them to the
multidimensional case is routine. The Jacobi identity imposes
conditions on all terms in the deformed bracket, which must be
satisfied degree by degree.

Thus, the search for the most general multiparameter deformation of
a given Lie algebra usually begins with the determination of the
space $H^2(\mathfrak{g};\mathfrak{g})$. An explicit basis given by
$2$-cocycles (representing the classes) determines an infinitesimal
deformation. One then tries to prolong this infinitesimal
deformation to all degrees. This prolongation method brings up the
{\it Massey brackets} which we will now describe, see \cite{Fu, FL,
Fia, Mill}.

Let $\mathfrak{g}_t$ be a $1$-parameter deformation of a Lie algebra
$\mathfrak{g}$, given by an infinitesimal cocycle $c=c^1$ and higher
degree terms $c^2$, $c^3$, \ldots . The Jacobi identity modulo
$t^{n+1}$ reads
\begin{equation} \label{br}
 \mathop{\sum}\limits_{i+j=n,\; i, j\geq 0}(c^i(c^j(x,y),z)+c^j(c^i(x,y),z)+{\rm
 cyclic}(x,y,z))=0,
\end{equation}
where ${\rm cyclic}(x,y,z)$ denotes the sum of all cyclic
permutations of the arguments of the expression written on the
left of it.

The expression \eqref{br} can be rewritten as
\[
 \mathop{\sum}\limits_{0\leq i,j\leq n;\;\; i+j=n}[[c^i,c^j]](x,y,z)=0,
\]
where the brackets $[[\cdot,\cdot]]$ are called  {\it Nijenhuis
brackets} (in differential geometry) or {\it Massey brackets} (in
deformation theory). The collection of the brackets
$[[\cdot,\cdot]]$ defines a graded Lie superalgebra structure on
$H^*(\mathfrak{g}; \mathfrak{g})$ (for examples, see \cite{LLS,
GL1}). The whole sum may then be expressed as a {\it Maurer\cherta
Cartan} equation:
\[
 \frac{1}{2}\sum_{i+j=n,\; i, j>0}[[c^i,c^j]]=dc^n
\]
because the term $[[c^0,c^n]]$ is just the left hand side of the
$2$-cocycle condition on $c^n$ in the Lie algebra $\mathfrak{g}$
(with adjoint coefficients) for the cochain $c^n$.

This gives a clear procedure for the prolongation of an
infinitesimal deformation (expressed here for simplicity only for a
$1$-parameter deformation): given a first degree deformation via a
cocycle $c=c^1$, one must compute its {\it Massey square} $[[c,c]]$.
If $[[c,c]]=0$, the infinitesimal deformation fulfills the Jacobi
identity and is thus a true deformation. If $[[c,c]]\in
Z^3(\mathfrak{g},\mathfrak{g})$ is not a coboundary, the
infinitesimal deformation is obstructed and cannot be prolonged. If
 $[[c,c]]=d\alpha$ with $\alpha\not=0$, then
$-\alpha t^2$ is the second degree term of the deformation. In order
to prolong to the third degree, one has to compute the next step ---
the Massey product $[[c,\alpha]]$. Once again, there are the three
possibilities $[[c,\alpha]]=0$, $[[c,\alpha]]=d\beta$ with
$\beta\not=0$ or $[[c,\alpha]]\not=d\beta$ for any $\beta$. If
$[[c,\alpha]]=d\beta$, then $\beta$ gives the third degree
prolongation of the deformation. In order to go up to degree $4$
then, one has to be able to compensate
$[[\alpha,\alpha]]+[[c,\beta]]$ by a coboundary $d\gamma$, and so
on. One must be careful to keep track of all terms coming in to
compensate low degree Massey brackets in a multiparameter
deformation.

The main difficulty in this kind of obstruction calculus is that the
representatives of the cohomology classes and the $\alpha$-,
$\beta$-, etc. cochains are not uniquely\footnote{If $\alpha$ is a
solution to the equation $d\alpha=[[c,c]]$, then $\alpha+
\mathrm{cocycle}$ is also a solution. The choice of a certain
$\alpha$ effects the expression of the $\beta$'s. The problem is how
to find a \lq\lq nice'' $\alpha$ in order to have as few
$\beta$-terms as possible and, more importantly, vanishing Massey
products in degrees higher than that of $\beta$. If we fail to
achieve this with $\alpha$, let us try to perform this with
$\beta$'s, and so on.} defined. A good choice of cochains may
considerably facilitate computations.

We computed cohomology and Massey products using Grozman's
\texttt{Mathematica}-based package \texttt{SuperLie}. The formula of
the following lemma was helpful in the computations. For any finite
dimensional Lie algebra $\mathfrak{g}$, all cochains with adjoint
coefficients may be expressed as sums of tensor products of the form
$x\otimes\omega$, where $x\in\mathfrak{g}$ and
$\omega\in\bigwedge^{\bcdot}(\mathfrak{g}^*)$. We are working with a
fixed basis of $\mathfrak{g}$ and the dual basis of
$\mathfrak{g}^*$.

\begin{Lemma}[Grozman] \label{lem1}
For any $c=a\otimes\omega$, where $x\in\mathfrak{g}$ and
$\omega\in\bigwedge^{\bcdot}(\mathfrak{g}^*)$, let $dc$ denote the
coboundary of $c$ in the complex with adjoint coefficients, while
$d\omega$ denotes the coboundary in the complex with trivial
coefficients and $da$ denotes the coboundary of $a\in\mathfrak{g}$
considered as a $0$-cochain in the complex with adjoint
coefficients. If $c=a\otimes\omega$, then $dc=a\otimes d\omega +
da\wedge\omega$.
\end{Lemma}

\begin{proof}
For any $x_1,\ldots,x_{p+1}\in\mathfrak{g}$, we have:
\[
\begin{split}
 dc(x_1,\ldots,x_{p+1})&=\sum_{1\leq i<j\leq p+1}(-1)^{i+j-1}a\otimes
 \omega([x_i,x_j],x_1,\ldots,\hat{x_i},\ldots,\hat{x_j},\ldots,x_{p+1})+{}\\
 &+\sum_{1\leq i\leq p+1}(-1)^{i}[x_i,a]\otimes
 \omega(x_1,\ldots,\hat{x_i},\ldots,x_{p+1})=\\
 &=(a\otimes d\omega)(x_1,\ldots,x_{p+1})+(da\wedge\omega)
 (x_1,\ldots,x_{p+1}).\qed
\end{split}
\]
\noqed\end{proof}

\subsection{In characteristic $2$} Let now $k$ be an algebraically closed field of
characteristic $2$. In this subsection $-1=1$, of course; the signs
are kept to make expressions look like in characteristics 0.

The vector space $\mathfrak{g}$ is a {\it Lie algebra} if endowed
with a bilinear map $[\cdot ,
\cdot]\colon\mathfrak{g}\times\mathfrak{g}\to\mathfrak{g}$
satisfying the Jacobi identity and anti-symmetry which for $p=2$
means $[x,x]=0$ for any $x\in\mathfrak{g}$. For vector spaces, the
wedge product is defined without a normalization factor:
\[
 a\wedge b=a \otimes b -b \otimes a.
\]

Grozman communicated to us the following definition of Lie algebra
cohomology in $\Char=2$ implemented in his \texttt{SuperLie}:

For $1$-cochains with trivial coefficients, the codifferential  is
defined as an operation dual to the Lie bracket:
\[
 d\colon \mathfrak{g}^*\rightarrow \mathfrak{g}^* \wedge \mathfrak{g}^*.
\]
For $q$-cochains with trivial coefficients, $d$ is defined via the
Leibniz rule. For cochains with coefficients in a module $M$, we
set
\[
\begin{split}
 &d(m):=-\mathop{\sum}\limits_{1\leq i\leq \dim \mathfrak{g}} [g_i, m] \otimes g_i^*,\\
 &d(m \otimes \omega):=d(m)\wedge\omega +  m \otimes d(\omega)
\end{split}
\]
for any $m\in M$, any $q$-cochain $\omega$, where $q>0$, and any
basis $g_i$ of $\mathfrak{g}$, cf. Lemma \ref{lem1}.

The Massey product is defined as follows:
\[
{}[a,b](x,y,z):= a(b(x,y),z)+b(a(x,y),z)+ \mathrm{ cyclic } (x,y,z),
\]
if $a$ and $b$ are non-proportional, whereas
\[
{}[a,a](x,y,z):=a(a(x,y),z)+ \mathrm{ cyclic } (x,y,z).
\]

\section{Main results: $p=3$}\texttt{}

\ssec{Deformations of $\fo(5)$: The results known}\texttt{} 1) A
{\bf 1-parameter} family of deformations of $\fo(5)$ is given by
Cartan matrices of $L(\eps)$.

Denote by $x$'s and $y$'s the Chevalley generators of $\fo(5)$.

 \sssbegin{Proposition} The Lie algebra $L(\eps)$ can be obtained as a deformation of
$\mathfrak{o}(5)$ generated by the $2$-cocycle $c_0$ bellow. The
bracket is as follows:
\[
{}[\cdot,\cdot]_{-1-\eps}=[\cdot,\cdot]-(1+\eps)\, c_0 + (1+\eps)^2
\alpha_0,
\]
where
\begin{equation}\label{czero}
\begin{array}{rcl} c_0&=& h_1\otimes \left(d x_2 \wedge d y_2\right)
+2\,
   h_1\otimes \left(d x_3 \wedge d y_3\right) +2\,
   h_1\otimes \left(d x_4 \wedge d y_4 \right) +
   h_2\otimes \left(d x_4 \wedge d y_4\right)\\
   && +
   x_1\otimes \left(d h_2 \wedge d x_1 \right) +2\,
   x_2\otimes \left(d h_2 \wedge d x_2\right) +
   x_2\otimes \left(d x_4 \wedge d y_3\right) +
   x_3\otimes \left(d x_4 \wedge d y_2 \right) \\
   &&+2\,
   x_4\otimes \left(d h_2 \wedge d x_4\right) +2\,
   y_1\otimes \left(d h_2 \wedge d y_1\right) +
   y_2\otimes \left(d h_2 \wedge d y_2 \right) +2\,
   y_2\otimes \left(d x_3 \wedge d y_4 \right)\\
   && +2\,
   y_3\otimes \left(d x_2 \wedge d y_4\right) +
   y_4\otimes \left(d h_2 \wedge d y_4 \right);\\[3mm]
\alpha_0 &=& h_1(dx_4 \wedge y_4).
\end{array}
\end{equation}
\end{Proposition}

\begin{proof}
Let us denote by $X$'s and $Y$'s the Chevalley generators of
$L(\varepsilon)$. The isomorphism is given by
\[
\renewcommand{\arraystretch}{1.4}\arraycolsep=2pt
\begin{array}{llll}
 X_i \longleftrightarrow x_i,
 & Y_i \longleftrightarrow y_i,
 & H_1 \longleftrightarrow h_1,
 & H_2 \longleftrightarrow h_2+(2-\varepsilon) h_1.
\end{array}
\]
\end{proof}

2) The {\bf 3-parameter} family of deformations of $\fo(5)$, denoted
here by $L(\eps, \delta,\rho)$, was constructed by Kostrikin (see
\cite{Kos}) as follows. Consider the contact Lie algebra
$\mathfrak{k}(3; \un)$, where $\un=(N_1,N_2,N_3)\in \Zee^3$,
generated by indeterminates $x$, $y$ and $t$ forming the algebra of
divided powers. As a vector space, $\mathfrak{k}(3; \un)$ is the
subspace of $\Kee[x,y,t]$ spanned by the monomials $x^iy^jt^k$ with
$0\leq i<p^{N_1}$, $0\leq j<p^{N_2}$, and $0\leq k < p^{N_3}$. As
usual in the divided power algebra, one has
\[
 w^i\cdot w^j=\binom{i+j}{i}w^{i+j}, \quad \text{and} \quad \partial_ww^i=w^{i-1}
\]
for $w=x,y$ or $t$. The contact bracket of polynomials $f$ and $g$
is defined by
\begin{equation}\label{cb}
 [f,g]=\triangle f\cdot\partial_t g - \partial_t f\cdot\triangle g
 +\partial_x f\cdot\partial_y g - \partial_y f\cdot\partial_x g
\end{equation}
with $\triangle f=2f - x\partial_x f - y\partial_y f$.

The {\it standard} $\Zee$-grading $\deg_{Lie}$ of $\mathfrak{k}(3;
\un)$ is defined by setting $\deg_{Lie}(f)=\deg f-2$, where
$\deg(x)=\deg(y)=1$ and $\deg(t)=2$. Then a basis of $L(\eps,
\delta,\rho)$ is given as follows:
\[
\renewcommand{\arraystretch}{1.4}
\begin{tabular}{|c|l|} \hline
$\deg$&the generator with weight $=$ its generating function \\
\hline \hline
$-2$&$E_{-2\alpha -\beta}=[E_{-\alpha}, E_{-\alpha-\beta}]=1;$\\
\hline
$-1$& $E_{-\alpha}=x;\quad E_{-\alpha-\beta}=[E_{-\beta}, E_{-\alpha}]=y;$   \\
\hline $0$& $H_{\alpha}=2 \eps t+ xy;\quad H_{\beta}=-xy; \quad
E_{\beta}=x^2; \quad
E_{-\beta}=-y^2;$   \\
\hline $1$& $E_{\alpha}=-(1+\eps)xy^2+ \eps
yt;\quad E_{\alpha+\beta}=[E_\beta, E_\alpha]=(1+\eps)x^2y+\eps x t;$   \\
\hline $2$& $E_{2\alpha+\beta}=[E_\alpha,
E_{\alpha+\beta}]=\eps(1+\eps)x^2y^2+ \eps^2
t^2.$   \\
\hline
\end{tabular}
\]
\normalsize The brackets involving new parameters are as follows
\begin{equation}\label{croc1}
\renewcommand{\arraystretch}{1.4}\arraycolsep=2pt
\begin{array}{lll}
 [E_{-2\alpha-\beta}, E_{-\alpha-\beta}]= \delta E_{\beta},\quad{}
 &[E_{-2\alpha-\beta}, E_{-\alpha}]= \rho E_{-\beta},
 &[E_{-2\alpha-\beta}, E_{-\beta}]= - \delta E_{\alpha+\beta},\\
 {}[E_{-2\alpha-\beta}, E_\beta] = \rho E_{\alpha},
 &[E_{-\alpha-\beta},E_{-\beta}]= -\frac{\delta}{\eps}
 E_{2\alpha+\beta},\quad{}
 &[E_{-\alpha}, E_{\beta}] = - \frac{\rho}{\eps} E_{2 \alpha+\beta}.
\end{array}
\end{equation}
\ssbegin{Remark} Kostrikin and Kuznetsov in \cite{KK} write
$H_\alpha$ as $t+xy$ instead.
\end{Remark}

\ssbegin{Proposition}\label{thm2} The Lie algebra $L(\eps,
\delta,\rho)$ can be obtained as a deformation of $\fo(5)$ generated
by the cocycle $(\ref{czero})$ and the following cocycles
\[
\begin{array}{ccl}
c_3&=& x_2\otimes \left(x_1^* \wedge y_4^*\right) +
   x_4\otimes \left(x_1^* \wedge y_2^* \right) +
   y_1\otimes \left(y_2^* \wedge y_4^* \right), \\
c_6&=&  x_1\otimes \left(y_3^* \wedge y_4^*\right) + 2 \,
   x_3\otimes \left(y_1^* \wedge y_4^* \right) +
   x_4\otimes \left(y_1^* \wedge y_3^* \right).
   \end{array}
\]
The bracket is as follows:
\[
\begin{array}{ccl}
[\cdot,\cdot]_t&=&[\cdot,\cdot]+t_1^2\, c_0 +t_3\, c_3 +t_4\, c_6 +
t_1^2 \alpha_0 + t_1t_3\, \alpha_3 + t_1t_4\, \alpha_6 +t_1^2t_3\,
\beta_3 + t_1^2t_4\, \beta_6+ \\
&&t_1^3t_4 \, \gamma + t_1^4t_4 \, \theta,
\end{array}
\]
where
\[
\eps=2-t_1, \quad \rho=\eps(\eps+2)\; t_3, \quad \delta= \varepsilon
(\varepsilon+2)(2+2 \varepsilon+\varepsilon^2)\; t_4.
\]
and
\[
\begin{split}
\alpha_6=- x_4\otimes \left(y_1^*\wedge y_3^* \right), \quad
   \alpha_3=-  x_4\otimes \left(x_1^* \wedge y_2^* \right),\quad
   \beta_3=-  x_2\otimes \left(x_1^*\wedge y_4^* \right) -y_1\otimes \left(y_2^*\wedge
   y_4^* \right), \\
   \beta_6=x_4\otimes \left(y_1^* \wedge y_3^*\right),\quad
   \gamma=- x_4\otimes \left(y_1^*\wedge y_3^* \right), \quad
   \theta= x_3\otimes \left(y_1^* \wedge y_4^*\right) -
    x_1\otimes \left(y_3^* \wedge y_4^* \right),
   \end{split}
\]
\end{Proposition}

\begin{proof} We can, of course, write down the whole multiplication
table of $\fo(5)_t$ but to make the paper shorter we will not do it.
Let us write only those constant structures for which we can deduce
the values of $\rho$ and $\delta$. Indeed
\begin{equation}  \label{cro2}
 \begin{split}
 &[y_2,x_2] =(2-t_1)x_3,\quad [y_2,x_1] = (t_1 t_3-t_3)x_4,\quad
 [y_4, y_3] = (t_1^4t_4+t_4)x_1.
\end{split}
\end{equation}

Since $L(\eps, \delta,\rho)$ was constructed in terms of $\fo(5)$,
then, deforming the bracket, it is natural to define the
isomorphisms between $L(\eps, \delta,\rho)$ and $\fo(5)_t$ as
follows:
\[
\renewcommand{\arraystretch}{1.4}\arraycolsep=2pt
\begin{array}{lllll}
 x_1 \leftrightarrow E_{\beta},
 & y_1 \leftrightarrow E_{-\beta},
 & h_1 \leftrightarrow H_{\beta},
 & x_2 \leftrightarrow E_{\alpha},
 & y_2 \leftrightarrow E_{-\alpha},\\
 h_2+(2-\eps) h_1 \leftrightarrow H_{\alpha},\quad{}
 & x_3 \leftrightarrow E_{\alpha+\beta},\quad{}
 & y_3 \leftrightarrow E_{-\alpha-\beta},\quad{}
 & x_4 \leftrightarrow E_{2\alpha+\beta},\quad{}
 & y_4 \leftrightarrow E_{-2\alpha-\beta}.\qed
\end{array}
\]
\noqed\end{proof}

3) Rudakov (cited in \cite{Kos}) constructed a {\bf 3-parameter}
family of deformations of $\mathfrak{o}(5)$ as the Cartan prolong of
the pair $(\Tee(a,b,c), \fgl(2))$. By construction, these deforms
linearly depend on parameters. 

\ssec{The deforms of $\fo(5)$: General picture}\texttt{} Since $\dim
H^2(\fo(5); \fo(5)))=5$, we will be dealing with five parameters,
denoted by $t_1,\ldots t_5$. We denote the Chevalley generators
corresponding to positive (resp. negative) roots by $x$ (resp. $y$).
The Lie algebra $\fo(5)$ has infinitesimal deformations given by the
following cocycles whose index is equal to their degree induced by
the $\Zee$-grading of $\fo(5)$ for which $\deg x_1=\deg x_2=1$ (here
$x_3=[x_1,x_2]$, $x_4=[x_2,x_3]$ and similarly for the $y$'s):
\[
\renewcommand{\arraystretch}{1.4}\arraycolsep=2pt
\begin{array}{ll}
 c_6 & =  x_1\otimes (y_3^* \wedge y_4^*)
   + 2  x_3\otimes (y_1^* \wedge y_4^*)+x_4\otimes (y_1^*\wedge y_3^*),\\
 c_{3} & =  x_2\otimes (x_1^*\wedge y_4^*)+
   x_4\otimes (x_1^*\wedge y_2^*)+
   y_1\otimes (y_2^*\wedge y_4^*),\\
   c_0 & = 2 h_1\otimes(x_2^* \wedge y_2^*)+
   2 h_1\otimes (x_3^* \wedge y_3^*)+
   2 x_1\otimes (x_3^* \wedge y_2^*)+
   y_1\otimes (x_2^* \wedge y_3^*),\\
 c_{-3} & = 2  x_1\otimes (x_2^*\wedge x_4^*)
   +  y_2\otimes (x_4^* \wedge y_1^*)+
   y_4\otimes (x_2^*\wedge y_1^*),\\
 c_{-6} &=   y_1\otimes (x_3^* \wedge x_4^*)
   +2 y_3\otimes (x_1^*\wedge x_4^*)+
   y_4\otimes (x_1^*\wedge x_3^*).
\end{array}
\]

Observe a symmetry between $c_6$ and $c_{-6}$, and between $c_3$ and
$c_{-3}$: there is an involution on the Lie algebra interchanging
$x$-generators and $y$-generators. One has to pay attention that
there is a sign involved $(2=-1)$ when passing from $c_3$ to
$c_{-3}$.

\ssbegin{Theorem} \label{thm1} The Lie algebra $\fo(5)$ admits a
$5$-parameter family of deforms denoted by $\fo(5,t)$, where
$t=(t_1, t_2, t_3, t_4, t_5)$.

The deformed bracket is defined by
\[
\begin{split}
 [\cdot,\cdot]_{t_1, t_2, t_3, t_4, t_5}&=
 [\cdot,\cdot] + t_1 \, c_0 + t_2\, c_{-3} + t_3 \, c_3 + t_4 \, c_6 + t_5\, c_{-6} +
 t_1 t_4 \, \alpha_{0,6} + t_1 t_2 \, \alpha_{0,-3}+{}\\
 & {}+ t_1 t_5 \, \alpha_{0,-6} + t_1 t_3 \, \alpha_{0,3} + t_4 t_5 \, \alpha_{6,-6} +
 t_2 t_3 \, \alpha_{-3,3} + t_1 t_4 t_5 \, \beta_{-6,0,6}+ t_1 t_2 t_3\, \beta_{3,0,-3},
\end{split}
\]
where
\[
\renewcommand{\arraystretch}{1.4}\arraycolsep=2pt
\begin{array}{llll}
 \alpha_{0,-6} & = y_1\otimes (x_3^*\wedge x_4^*),&
   \alpha_{0,6} & = x_1\otimes (y_3^* \wedge y_4^*),\\
   \alpha_{0,3} & =  y_1\otimes (y_2^*\wedge y_4^*),\quad{}&
   \alpha_{0,-3} &= 2 x_1\otimes (x_2^* \wedge x_4^*),\\
\end{array}
\]
\[
\tiny
\renewcommand{\arraystretch}{1.4}\arraycolsep=2pt
 \begin{array}{ll}
 \alpha_{-3,3}  = & 2 h_2\otimes (x_2^* \wedge y_2^*) +
     h_2\otimes (x_3^* \wedge y_3^*) -
     h_2\otimes (x_4^* \wedge y_4^*)+
    2 x_1\otimes (x_3^* \wedge y_2^*) +
      x_2\otimes (h_1^* \wedge x_2^*)+{}\\
   &{}+ 2 x_2\otimes (x_3^* \wedge y_1^*) +
    x_3 \otimes (h_1^* \wedge x_3^*) +
   2 x_4 \otimes (h_1^* \wedge x_4^*) +
    x_4\otimes (x_2^* \wedge x_3^*) +
  2 y_2\otimes (h_1^* \wedge y_2^*)+{}\\
   &{}+ 2 y_2\otimes (x_3^* \wedge y_4^*)+
   2 y_3\otimes (h_1^* \wedge y_3^*) +
   2 y_3\otimes (x_2^*  \wedge y_4^*)-
   2 y_3\otimes (y_1^*  \wedge y_2^*) -
   2 y_4\otimes (h_1^* \wedge y_4^*),\\[3mm]
 \alpha_{-6,6} = &   h_2\otimes (x_2^* \wedge y_2^*) +
   2 h_2\otimes (x_3^* \wedge y_3^*) +
    h_2\otimes (x_4^* \wedge y_4^*) +
   2 x_1\otimes (x_3^* \wedge y_2^*)+
   2x_2\otimes (h_1^* \wedge x_2^*)+{} \\
   &{}+ 2 x_2\otimes (x_3^* \wedge y_1^*) +
   2 x_3\otimes (h_1^* \wedge x_3^*) +
    x_4\otimes (h_1^* \wedge x_4^*) +
    x_4\otimes (x_2^* \wedge x_3^*) +
    y_2\otimes (h_1^* \wedge y_2^*)+{}\\
   &{}+ 2  y_2\otimes (x_3^* \wedge y_4^*) +
    y_3 \otimes (h_1^* \wedge y_3^*) +
   2 y_3\otimes (x_2^* \wedge y_4^*) +
    y_3\otimes (y_1^* \wedge y_2^*) +
   2  y_4\otimes (h_1^* \wedge y_4^*),
\end{array}
\]
and
\[
\tiny
\renewcommand{\arraystretch}{1.4}\arraycolsep=2pt
 \begin{array}{ll}
 \beta_{-6,0,6} = &2 h_1\otimes (x_2^* \wedge y_2^*) +
    h_1\otimes (x_3^*   \wedge y_3^*) +
   2 h_1\otimes (x_4^* \wedge y_4^*) +
   2 x_1\otimes (x_3^* \wedge y_2^*) +
   2 x_2\otimes (h_2^* \wedge x_2^*)+{}\\
   &{} +2 x_3\otimes (h_2^* \wedge x_3^*) +
   x_4\otimes (h_2^*  \wedge x_4^*) +
    y_2\otimes (h_2^* \wedge y_2^*) +
    y_3\otimes (h_2^* \wedge y_3^*) +
   2 y_4\otimes (h_2^* \wedge y_4^*),\\[3mm]
 \beta_{-3,0,3} = &  h_1\otimes (x_2^* \wedge y_2^*) +
   2 h_1\otimes (x_3^* \wedge y_3^*) +
   h_1\otimes (x_4^* \wedge y_4^*) +
   2 x_1\otimes (x_3^* \wedge y_2^*) +
   x_2\otimes (h_2^* \wedge x_2^*)+{}\\
   &{} + x_3\otimes (h_2^* \wedge x_3^*) +
   2 x_4\otimes (h_2^* \wedge x_4^*) +
   2 x_4\otimes (x_2^* \wedge x_3^*) +
   2 y_2\otimes (h_2^* \wedge y_2^*) +
    y_2\otimes (x_3^* \wedge y_4^*)+{}\\
   &{} +2 y_3\otimes (h_2^* \wedge y_3^*) +
    y_3\otimes (x_2^* \wedge y_4^*) +
   y_4\otimes (h_2^* \wedge y_4^*).
\end{array}
\]
\end{Theorem}

\begin{proof} The proof is a direct computation assisted by Grozman's
\verb"Mathematica"-based package \verb"SuperLie" (\cite{Gr}). We
compute the Massey brackets in each degree and check if this bracket
is a coboundary. For example, we can easily get
\[
\begin{split}
 &[[c_{-3}, c_{-6}]]= 0,\qquad [[ c_{3}, c_{3}]]= 0,\\
 &[[c_{0}, c_{6}]] =  2  h_1\otimes (y_1^* \wedge y_3^*\wedge y_4^*) +
2 x_1\otimes (y_1^* \wedge y_2^* \wedge y_4^*) + 2  x_3\otimes
(x_2^* \wedge y_3^* \wedge y_4^*).
\end{split}
\]
Besides, we can show that $[[ c_{0}, c_{6}]]=-d\alpha_{0,6} $, where
$\alpha_{0,6}$ is as above. It is here that we use the formula of
Lemma \ref{lem1} in order to compute $d\alpha_{0,6}$. Indeed, we
have:
\[
\renewcommand{\arraystretch}{1.4}\arraycolsep=2pt
\begin{array}{rl}
 d\alpha_{0,6}= & d( x_1\otimes y_3^* \wedge y_4^*)= \\
  =&   dx_1\otimes y_3^* \wedge y_4^*  +
       x_1\otimes dy_3^* \wedge y_4^* +
     2 x_1\otimes y_3^* \wedge dy_4^*= \\
  =&  x_3\otimes x_2^*\wedge y_3^*\wedge y_4^* +
       h_1\otimes y_1^*\wedge y_3^*\wedge y_4^* +
      x_1\otimes h_1^*\wedge y_3^*\wedge y_4^* +{} \\
  &{} + 2 x_1\otimes h_2^*\wedge y_3^*\wedge y_4^* +
       x_1\otimes y_1^*\wedge y_2^*\wedge y_4^* +  x_1\otimes y_3^*\wedge h_1^*\wedge y_4^*
       +{}\\
  &{}  + x_1\otimes y_3^*\wedge y_4^*\wedge h_2^* = \\
  =&  x_3\otimes x_2^*\wedge y_3^*\wedge y_4^* +
       h_1\otimes y_1^*\wedge y_3^*\wedge y_4^* +
       x_1\otimes y_1^*\wedge y_2^*\wedge y_4^*.
\end{array}
\]
In this computation we used the knowledge of the explicit form of
$dx_1$ (the coboundary of a $0$-cochain with adjoint coefficients)
which we extracted from the multiplication table:
\[
 dx_1=x_3\otimes x_2^* + h_1\otimes y_1^* + y_2\otimes y_3^* +
 x_1\otimes h_1^* + 2 x_1\otimes h_2^*,
\]
and the explicit form of $dy_3^*$ and $dy_4^*$ (the coboundary of
a $1$-cochain with values in the ground field):
\[
\begin{split}
 &dy_3^*=y_1^*\wedge y_2^*  + y_4^*\wedge x_2^* + y_3^*\wedge h_1^*,\\
 &dy_4^*=2y_4^*\wedge h_2^* + 2 y_3^*\wedge y_2^*.
\end{split}
\]
We have chosen the cocycles so that their Massey squares are $0$. As
explained above, the $\alpha$'s are not unique. We hoped that we can
choose them so that the $\beta$'s (corresponding to Massey products
of degree three) are ALL zero. Unfortunately, this is not possible.
Nevertheless, we can choose the $\alpha$'s so that a large number of
the $\beta$'s vanish. Once this is done, we can deal with the
$\beta$'s. Rather long computations with the remaining free
parameters in degree $4$ show that the $\alpha$- and
$\beta$-cochains can be chosen so that all Massey brackets in degree
$4$ vanish. This was our choice.
\end{proof}

\section{Main results: $p=2$}\texttt{}

\ssec{Deforms of $\fo^{(1)}(5)$}\texttt{} In this subsection, $p=2$,
and hence the orthogonal Lie algebra $\fo(5)$ is not simple and of
dimension 15. The Lie algebra $\fo^{(1)}(5)$, the derived of
$\fo(5)$, is simple and of dimension 10. It is realized by means of
the Cartan matrix and the generators (same with the $y$'s):
\begin{equation*}
 \mat{
 1 & -1 \\
 -1 & 0
 };
 \quad  x_1, x_2, x_3=[x_1,x_2],\ \ x_4=[x_1, x_3];
\end{equation*}
see \cite{BGL}. From \cite{BGL4} we know that $\dim
H^2(\fo^{(1)}(5); \fo^{(1)}(5))=4$, so we will be dealing with four
parameters, denoted by $t_1,\ldots, t_4$.

The Lie algebra $\mathfrak{o}^{(1)}(5) $ has infinitesimal
deformations given by the following cocycles:
\[
\tiny
\renewcommand{\arraystretch}{1.4}\arraycolsep=2pt
 \begin{array}{lcl}
 c_{4} &= & h_1\otimes (y_2^*\wedge y_4^*) + x_1\otimes (y_2^*\wedge y_3^*)+
    x_2\otimes (h_2^*\wedge y_4^*) + x_3\otimes(h_2^*\wedge y_3^*) +
    x_4\otimes (h_2^*\wedge y_2^*) + y_1\otimes(y_3^*\wedge y_4^*),\\[3mm]
 c_{-2} & = & h_1\otimes ( x_4^* \wedge y_2^* ) + x_2\otimes(h_1^*\wedge x_4^* ) +
    x_2\otimes(h_2^* \wedge x_4^* ) + x_3\otimes (x_1^* \wedge x_4^* ) +
    y_1\otimes(h_1^* \wedge x_1^* ) + y_1\otimes (h_2^* \wedge x_1^* ) +{}\\
   && {}+ y_1\otimes (x_3^* \wedge y_2^* ) + y_4\otimes (h_1^* \wedge y_2^* )
   + y_4\otimes(h_2^* \wedge y_2^* ) + y_4\otimes(x_1^* \wedge y_3^* ),\\[3mm]
 c_{2} &= & h_1\otimes (x_2^* \wedge y_4^* ) + x_1\otimes(h_1^* \wedge y_1^* ) +
    x_1\otimes (h_2^* \wedge y_1^* ) + x_1\otimes (x_3^* \wedge y_4^*) +
    x_3\otimes (x_2^* \wedge y_1^* ) + x_4\otimes (h_1^* \wedge x_2^*) +{}\\
  &&{}+ x_4\otimes (h_2^* \wedge x_2^* ) + y_2\otimes (h_1^* \wedge y_4^* ) +
  y_2\otimes (h_2^* \wedge y_4^* ) + y_2\otimes (y_1^* \wedge y_3^* ),
\end{array}
\]
and the 2-cocycle $c_{-4}$ is obtained from $c_{4}$ by changing $x$
by $y$ and $y$ by $x$. These cocycles $c_i$ are chosen so that
$[[c_i, c_i]]=0$ (which, fortunately, is possible) and having
shortest possible expression (for esthetic reasons).

\ssbegin{Theorem}\label{thm3} The Lie algebra $\fo^{(1)}(5)$ admits
a $4$-parameter family of deformations denoted by $\fo^{(1)}(5; t)$,
where $t=(t_1,t_2,t_3,t_4)$. The deformed bracket is given by the
formula
\[
\begin{split}
 [\cdot ,\cdot]_{t_1,t_2, t_3, t_4} ={} &
    [\cdot,\cdot]+t_1 c_{-4} + t_2 c_{4} + t_3 c_{-2} + t_4  c_2 +
    t_1 t_3 \alpha_{-4,-2} + t_1 t_4 \alpha_{-4,2} + t_2 t_3 \alpha_{4,-2} +{}\\
 & + t_2 t_4 \alpha_{4,2} + t_3 t_4 \alpha_{-2,2} + t_2 t_3 t_4 \beta_{4,-2,2} +
    t_3^2\, t_4 \beta_{-2,-2,2} + t_3 t_4^2 \beta_{2,-2,2} +{}\\
 & + t_3^2 t_4^2 \alpha_{-2,2} + t_3^3 t_4^2 \beta_{-2,-2,2} +
    t_3^2 t_4^3 \beta_{2,-2,2} + t_3^3 t_4^3 \varrho,\\
\end{split}
\]
where
\[
\renewcommand{\arraystretch}{1.4}\arraycolsep=2pt
\begin{array}{lcllcl}
 \alpha_{-4,-2} & = & y_4\otimes(h_2^* \wedge x_4^*),
    &
    \alpha_{-4,2} & = & y_2\otimes(h_2^* \wedge x_2^* ),\\
 \alpha_{4,-2} & = & x_2\otimes (h_2^* \wedge y_2^* ) + y_1\otimes(y_2^*\wedge y_3^*),\quad{}
    &
    \alpha_{4,2} & = & x_1\otimes(y_3^* \wedge y_4^* ) + x_4\otimes(h_2^* \wedge y_4^* ),\\
 \beta_{2,-2,2} & = & h_1\otimes x_2^* \wedge y_4^* ) + x_1\otimes (x_3^* \wedge y_4^*),
    &
    \beta_{-2,-2,2} & = & h_1\otimes(x_4^* \wedge y_2^* ) + y_1\otimes(x_3^* \wedge y_2^* ),\\
 \beta_{4,-2,2} & = & x_3\otimes(h_2^* \wedge y_3^* ),
    &
    \alpha_{-2,2} & = & h_1\otimes (x_2^*\wedge y_2^* ) + h_1\otimes(x_4\wedge y_4^*)+{}\\
    &&&&& {}+ x_1\otimes(x_3^* \wedge y_2^* )+y_1\otimes (x_3^* \wedge y_4^*),
\end{array}
\]
\[
\begin{split}
 \varrho ={} & x_1\otimes (x_4^* \wedge y_3^* ) + x_2\otimes (h_2^*\wedge x_2^* ) +
    x_3\otimes (h_2^* \wedge x_3^* ) + x_4\otimes (h_2^* \wedge x_4^* ) +
    y_1\otimes (x_2^* \wedge y_3^* ) +{}\\
    &{} +y_2\otimes (h_2^* \wedge y_2^* ) + y_3\otimes (h_2^* \wedge y_3^*) + y_4\otimes (h_2^* \wedge y_4^* ).
\end{split}
\]
\end{Theorem}
\begin{proof} We follow the proof of Theorem \ref{thm1} {\it mutatis
mutandis}.\end{proof}

\section{Remarks}
\ssec{On the Lie algebra $N$ of \cite{KK}} Kostrikin and Kuznetsov
 made a mistake, acknowledged by
M.~Kuznetsov, so one should disregard all about the Lie algebra $N$
in \cite{KK}.

\ssec{Warning: The fact that $H^2(\fg;\fg)\neq0$ is not equivalent
to deformability of $\fg$}\label{alapp}

Rudakov writes in \cite{Ru} (Corollary 1) that a given Lie algebra
has no non-trivial deformations if and only if $H^2(\fg;\fg)\neq 0$.
This is wrong.

There is no one-to-one correspondence between cohomology classes of
$H^2(\fg;\fg)$ and deformations: There can be more deforms than
$\dim H^2(\fg;\fg)$ (see \cite{FF}), in positive characteristic
there can be fewer deforms than $\dim H^2(\fg;\fg)$ --- even though
all infinitesimal deforms might happen to be integrable --- if the
structure constants are not rational, as is written in \cite{KK}.
Here we will illustrate this statement of a too short and intriguing
note \cite{KK} by a simple example.

Let $\Kee$ be a field of characteristic $p>0$ such that each of its
element has a root of degree $p$ (for example, any finite field).
Let $\fg$ be a Lie algebra over $\Kee$ with basis
$(e_0,\dots,e_{p-1},f)$ and the following commutation relations
(here $i+1$ in the subscript should be understood modulo $p$):
$$
{}[e_i,e_j]=0;\quad [f,e_i]=e_{i+1}\qquad \text{for
any~}i,j=0,\dots,p-1.
$$

Let us consider the following deformation of $\fg$ with parameter
$a$:
\begin{equation}\label{deform}
\begin{array}{lcll}
{}[e_i,e_j]_a&=&0&\text{for any~}i,j=0,\dots,p-1;\\
{}[f,e_i]_a&=&e_{i+1}& \text{for any~}i=0,\dots,p-2;\\
{}[f,e_{p-1}]_a&=&(1+a)e_0.&
\end{array}
\end{equation}

\begin{Claim} The deformation $(\ref{deform})$ of the initial bracket is
trivial.\end{Claim}

\begin{proof} Let us consider linear operators $A_a$ which act as
follows on basic elements of $\fg$:
\begin{equation}\label{op}
A_af=\sqrt[p]{1+a}f;\quad A_ae_i=(\sqrt[p]{1+a})^ie_i\qquad\text{for
all~}i=0,\dots,p-1.
\end{equation}
Then
$$
{}[A_ax,A_ay]=A_a[x,y]_a\qquad\text{for all~}x,y\in\fg.
$$
\end{proof}The 2-cocycle corresponding to the infinitesimal
version of the deformation (\ref{deform}) is proportional to
$$
z=e_0\otimes \psi\wedge\phi_{p-1},
$$
where
$(\phi_0,\dots,\phi_{p-1},\psi)$ is the basis of $\fg^*$ dual to
$(e_0,\dots,e_{p-1},f)$.

\begin{Claim} The element $z$ represents a nontrivial cocycle of
$H^2(\fg;\fg)$.\end{Claim}

\begin{proof} The algebra $\fg$ has a $\Zee/p\Zee$-grading such
that
$$
\deg f=1;\quad\deg e_i=i\qquad\text{for~}i=0,\dots,p-1.
$$
The degree of $z$ in the corresponding grading of $C^*(\fg;\fg)$ is
equal to $0$. The subspace of $C^1(\fg;\fg)$ of degree $0$ is
spanned by the elements
$$
\begin{array}{ll}
e_i\otimes \phi_i&\text{for }i=0,\dots,p-1;\\
 e_1\otimes\psi, \quad f\otimes\phi_1, \quad  f\otimes\psi.
 \end{array}
$$
So $B^2(\fg;\fg)_0$ is the linear span of the elements (here $i\pm
1$ or $j\pm 1$ in the subscript should be understood modulo $p$)
\begin{equation}\label{B20}\begin{array}{ll}
d(e_i\otimes\phi_i)=e_{i+1}\otimes
\phi^i\wedge\psi-e_i\otimes\phi_{i-1}\otimes\psi&\text{for~}i=0,\dots,p-1;\\
d(f\otimes\phi_1)=\sum\limits_{i\neq 1}e_{i+1}\otimes
\phi_i\wedge \phi_1-f\otimes\phi_0\wedge\psi;&\\
d(f\otimes\psi)=\sum\limits_{i=0}^{p-1}e_{i+1}\otimes
\phi_i\wedge\psi&.
\end{array}
\end{equation}
(The differential $d(e_1\otimes\psi)$ vanishes, so does not count.)
Let us consider the linear map $L:C^2(\fg;\fg)\to\Kee$ defined on
the basic 2-cochains (here $i,j,k=0,\dots,p-1$) as follows:
$$
\begin{array}{ll}
L(e_i\ot \phi_j\wg \phi_k)=0;&L(e_i\otimes \phi_j\wg\psi)=\delta_{i,j+1};\\
L(f\ot \phi_i\wg \phi_j)=0;&L(f\ot \phi_i\wg\psi)=0.
\end{array}
$$
The value of $L$ on all the elements (\ref{B20}) is equal to $0$,
but $L(z)=-1$. Thus $z$ is not a linear combination of the elements
(\ref{B20}).
 \end{proof}

If the deformation (\ref{deform}) is trivial, why can $z$ not be
obtained  as the differential of $1$-cochain $C$ such that
$C(x)=\pderf{A_a x}{a}$ (whatever such partial derivative might mean
in characteristic $>0$) for any $x\in\fg$? This is because in
characteristic $p$ the function $\sqrt[p]{1+a}$ is not
differentiable (again, whatever \lq\lq differentiable" means here).

The situation is opposite, in a way, to the one described in
\cite{FF} for infinite dimensional Lie algebras over algebraically
closed fields of characteristic 0.

\label{lastpage}

\end{document}